\documentclass{article}\usepackage{newcent}\usepackage{amsfonts}\usepackage{amsmath}\usepackage{graphicx}\begin{document}\large
\def \m {\vspace{6pt}

\noindent} 

\vfill\eject \centerline{Functorial affinization of Nash's manifold}
\m {\small  
\m\pagenumbering{roman}{\it Abstract.} Let $M$ be a singular irreducible complex manifold of
dimension $n.$ There are ${\mathbb Q}$ divisors $D[-1], D[0], D[1],...,D[n+1]$
on Nash's manifold $U\to M$ such that $D[n+1]$ is relatively ample on
bounded sets, $D[n]$ is relatively eventually basepoint free on
bounded sets, and $D[-1]$ is canonical with the same relative
plurigenera as a resolution of $M.$ The divisor $D=D[n]$ is the 
supremum of divisors ${1\over i}D_i.$ An arc containing one
singular point of $M$ lifts 
to $U$ if and only if the generating number of $\oplus_i {\cal O}_{\gamma}(D_i)$ is finite. When finite it equals  $1+(K_U-K)\cdot \gamma$ where ${\cal O}_U(K)$
is the pullback mod torsion of $\Lambda^n\Omega_M.$ If $C$ is a
 complete curve in $U$ then ${{-1}\over{n+1}}K_U\cdot C=D_1\cdot C + 
D_{n+2} \cdot C + D_{(n+2)^2}\cdot C +.....$ When there are infinitely 
many nonzero terms the sum should be taken formally or $p$-adically for
a prime divisor $p$ of $n+2.$   There are finitely many  nonzero terms
 if and only if $C\cdot D=0.$  
The natural holomorphic map $U\to M$ factorizes through the 
 contracting map
$U\to Y_0.$ The Grauert-Riemenschneider sheaf of $M,$ if $M$ is bounded, agrees with ${\cal H}om({\cal O}_M(D_{(n+2)^i-1}),\  {\cal O}_M(D_{(n+2)^i}))$ for large $i.$  
If $M$ is projective,  singular $s$-dimensional foliations
on $M$ such that $K+(s+1)H$ is a finitely-generated divisor of Iitaka dimension one are completely resolvable, where $K$ is the canonical divisor of the foliation. 
\m According to a question of [7] it is not known whether $Y_0$ always has canonical singularities.  
It is not known whether the analogue of Nash transforms, locally principalizing  the reflexivication of $\Lambda^n\Omega_M,$ 
eventually
converges, but this is true for toric surfaces. Conjecturally 
a  split faithful  action of a  commutative Lie algebra can never be
completely  resolved
if the weights do not form a basis (possibly with multiplicity) of the 
dual of the Lie algebra.  This is not known. It is not understood
under what conditions  $Y_0$ can be connected to the relative 
canonical model
by proper  maps over $M.$ 
\m \m\m\m\m\m\m\m\m\m \m \hfill $\begin{matrix}

 \hbox{\it John Atwell Moody}\hfill \cr 
\hbox{\it Coventry}\hfill \cr\hbox{\it April 12, 2010}\hfill \end{matrix}$}
\vfill\eject
\pagenumbering{arabic}
\m \centerline{0. Preliminaries}
\m {\bf 0.1  Hypotheses.} Throughout this article, $M$ will be a singular
complex manifold; that is, a  
Hausdorff space with a countable open
cover by closed analytic subspaces of complex domains, furnished
with the reduced structure sheaf.
 We suppose that $M$ is irreducible and let $n=dim(M).$ An open subset
of $M$ will be called {\it bounded} if its closure is compact.

\m {\bf 0.2 The Nash manifold.} Let $...M_2\to M_1 \to M_0=M$ be
the sequence of Nash blowups of $M$ [4] and let $U_0 \subset U_1 \subset U_2
\subset ...$ be the ascending chain of open immersions, where $U_i$
is the regular locus of $M_i.$ Let $U$ be the colimit, so $U$ is the
ascending union of the $U_i;$ every point of $U$ is contained in
one of the $U_i$ and the $U_i$ are complex manifolds; and so $U$ is
also a complex manifold. It is, and is here also defined to be,
Nash's manifold.  It is not known whether the structural map $U\to M$
is always proper, and this is a question raised by Nash's work. 
\m {\bf 0.3 Canonical Singularities.} 
For each natural number $r$ let $V_r \to M$ be the universal singular
manifold such that the pushforward of  $\omega^{\otimes r}$ from a
 resolution of $V_r$ is invertible and relatively basepoint free
over $M.$   
The $V_i$ have inclusions according to divisibility
$$\begin{matrix}
&&V_2&&\cr
&\nearrow&&\searrow \cr
V_1&&&&V_6 \cr
&\searrow&&\nearrow\cr
&&V_3&&\end{matrix}
$$
etc. Set $V$ to be the colimit of the $V_i.$ The maps $V_i\to V$
are open immersions providing
an open cover. Since each $V_i$ has
canonical singularities, $V$ itself is a manifold
with canonical singularities. That is to say,  if $s:V'\to V$ is a resolution then every point of $V$ has
an index $r$ and a Stein neighbourhood $T$
such that the $\Gamma(s^{-1}T,\omega^{\otimes r} )$
is a principal  $\Gamma(T, {\cal O}_V)$ module.
As in the case of the
Nash manifold,
there is a natural map $V\to M,$ conjectured to be a proper
map [7], now known when the singularities of $M$ are isolated,
and much more generally ([13,21,22,23] further discussion in the appendix).

\vfill\eject\noindent {\bf 0.4 Constructions of Grothendieck and Atiyah}
\m Let now ${\cal F}$ be any coherent analytic sheaf on $M,$ we
assume ${\cal F}$ is torsion-free and of rank one. 
Let $f:{\bf V}({\cal F})\to M$ be
Grothendieck's map of singular manifolds [2], definition 1.7.8, with analytic action
of $Gl_1({\mathbb C})$ freely and transitively on the fibers
with fixed locus $M$ itself, which is characterised by the property
that for $i\in {\mathbb Z},$ 
the coherent sheaf of complex valued holomorphic functions on    $F$
which transform according
to the 
character of degree $i$ are just $0$ if $i<0$ and the power ${\cal F}^{\otimes i}/torsion$ otherwise. 
If ${\cal F}$ is invertible, then ${\bf V}({\cal F})$ can be identified point-by-point
with the dual of the line bundle whose section sheaf is ${\cal F}.$ 
In general, [2] paragraph 1.7.9,  the section sheaf of ${\bf V}({\cal F})$ is ${\cal H}om({\cal F},
{\cal O}_M).$
The subsheaf of $f_*\Omega_{ {\bf V}({\cal F})}/torsion$ which transform according to the 
character of degree one is called the {\it principal parts sheaf 
of $ {\cal F} $ modulo torsion } and will be denoted ${\cal P}({\cal F}),$
and the exact sequence $0\to {\cal F}\otimes \Omega_M/torsion 
\to {\cal P}({\cal F})\to {\cal F}\to 0$ is Atiyah's sequence [1]. 
There is a natural splitting $i$  of Atiyah's sequence as a sheaf
of complex vector-spaces, and the sheaf of  ${\cal O}_M$-module splittings 
assigns to an open set $T \subset M$ precisely the 
 $i+\nabla$ where $\nabla:{\cal F}_{|T }\to {\cal P}({\cal F}_{| T})$ runs over the distinct  connections on the restriction of
${\cal F}$ to $T.$ 
\m \vfill\eject\centerline{1. The ordering on coherent sheaves} 
\m\m {\bf 1.  Theorem} (see [17])  If ${\cal F}$ and ${\cal G}$ are two such
coherent sheaves then there is a natural transformation
$${\cal G}^{n+1}\Lambda^{n+1}{\cal P}({\cal F})\to \Lambda^{n+1}{\cal P}({\cal GF})\ \ \ (1) $$
such that for each ${\cal F}$ and ${\cal G},$ 
the map  pulls back on $Bl_{\cal FG}(M)$ to the natural map (where $\tau:Bl_{{\cal G}{\cal F}}M\to Bl_{\cal F}M$ over $M$)
$$ \tau^*\Lambda^n\Omega_{Bl_{\cal F}M}/torsion \to \Lambda^n\Omega_{Bl_{\cal FG}M}/torsion$$
twisted by the $n+1$ power of the invertible sheaf which results
by pulling back ${\cal FG}$ and reducing mod torsion. The map is not
an $n+1$ exterior power.
It generalizes the 
special case when $Bl_{\cal F}M$ is normal.  
\m {\bf 2. Corollary.} Suppose $Bl_{\cal F}M$ is normal. Then, on each bounded
open subset of $M,$  ${\cal G}$ is a divisor of a power of ${\cal F}$
as sheaves of fractional ideals 
if and only if (1)/torsion becomes surjective 
 after multiplying by 
a power of ${\cal FG}.$
\m 
Among multiplicatively closed `sets' of torsion-free coherent sheaves on $M$ closed
under  division and multiplication, those which happen to be finitely generated
on bounded open sets of $M$ ,
are generated by a single element on each bounded open set. Such multiplicatively closed
sets of sheaves correspond bi-uniquely with singular manifolds
$N\to M$ over $M$ whose structure map is locally projective of degree one. The 
corollary for example can be interpreted as saying that on bounded
open sets  the map (1) is an isomorphism
 after multiplying by some power of ${\cal FG}$ if and only if ${\cal G}$ 
is already contained in the smallest multiplicatively and divisbilitively
closed set containing ${\cal F}.$ 
\m Stein factorization implies that finitely generated  multiplicatively and divisibilitively
closed sets of torsion free coherent sheaves of rank one have the 
descending chain property on bounded open subsets of $M.$ The unique
minimal set for example is the set of invertible sheaves.
The corresponding chain condition can be included in the
conclusion of the theorem.
\m \vfill\eject\centerline{2. A  basepoint free theorem}
\m  Let
$f:N\to M$ be a locally projective morphism of degree one. Generalizing
from the case of invertible sheaves, we say a torsion-free coherent sheaf
of rank one ${\cal F}$ on $N$ is {\it spanned relative to $f$} if each point
$p$ of $M$ has a neighbourhood such that the restriction of  
${\cal F}$ to the
inverse image of the neighourhood is generated by global sections; and we say 
that ${\cal F}$ is {\it very ample relative to $f$} if all 
${\cal F}^{\otimes i}/torsion $  are spanned
relative to $f$ and  the meromorphic map $N - \to  {\cal P}roj\oplus_i(f_*{\cal F})^i$ is
the inverse of a morphism. 
\m {\bf 3. Corollary.} 
Let ${\cal L}$ now be an invertible sheaf on $N$ which is very ample
relative to $f.$ Then ${\cal L}^{n+1}\Lambda^n\Omega_N/torsion$ is spanned
relative to $f$ and ${\cal L}^{n+2}\Lambda^n\Omega_N/torsion$ is very
ample relative to $f.$
\m The inverse morphism is the Nash blowup, which makes its appearance
in this  way.
\m
The primary difference between Miles' model of [7] and the Nash model
is that Miles' model has a functorial and canonical affinization.
This difference disappears upon contemplation of the formula
for the partial sum of a geometric series
$$(n+2)^s=(n+1)(1+(n+2)+...+(n+2)^{s-1}).$$
Namely, starting from a torsion-free coherent sheaf ${\cal F}$
on $M,$ we define for each natural number $i$ a new torsion
free coherent sheaf, functorial with respect to Grothendiecks' category 
of holomorphic maps of $M$ 
and coherent sheaf maps of ${\cal F},$
by the rule that if the base-$(n+2)$-expansion of $i$ is
$a_0+a_1(n+2)+...+a_s(n+2)^s$ with $0\le a_i<(n+2),$ 
then we write ${\cal F}_i={\cal F}_1^{\otimes a_0}\otimes...\otimes{\cal F}_{(n+2)^s}^{
\otimes a_s}/torsion$
and when $i=(n+2)^j$ is a power of $n+2$ we write 
${\cal F}_i=\Lambda^{n+1}{\cal P}({\cal F}\otimes {\cal F}_{{i-1}\over{n+1}})/torsion.$
The morphism of theorem 1. provides the needed `carrying' map
$${\cal F}_{(n+2)^j}^{\otimes (n+2)}\to {\cal F}_{(n+2)^{j+1}}$$
such that the ${\cal F}_i$ can be multiplied following the
usual manner in integer expressions to the base $n+2$ are
added together.
\m {\bf 4. Theorem.} If ${\cal F}$ is a very ample invertible sheaf
on $M$  then the ${\cal F}_i$ are torsion free coherent sheaves on $M$  which
are generated by global sections.
 \vfill\eject\centerline{3.  The next-to last stage of Nash's tower}
\m\m 
\m {\bf 5. Theorem.} 
Suppose the Nash tower  $...M_{m+1}\to M_m \to ... \to M_0=M$
satisfies $M_{m+1}=M_m.$ Then $U=M_m,$  $\oplus_i{\cal F}_i$ is finite type 
 (although
the local generating degree can be larger than $(n+2)^m$ and might not be 
bounded)
 and taking  $Y={\cal P}roj
 \oplus {\cal F}_i$  there is
a pullback diagram
$$\begin{matrix} M_m&\to &M_{m-1}\cr
\downarrow &&\downarrow\cr
Y&\to&M\end{matrix}.$$
In other words,  $U=M_m$ can always be built by pulling back the penultimate
term of the Nash tower $M_{m-1}\to M$ along $Y\to M.$ In general, 
by section 6, properness of $U\to Y$ is equivalent to $\oplus {\cal F}_i$
being (locally) of finite type. \m
Proof.  ${\cal F}_{(n+2)^m-1}
={\cal F}_{(n+2)^{m-1}-1}\otimes
{\cal F}_{(n+2)^{m-1}}^{n+1}/torsion$  while  $M_m, M_{m-1}, Y$ are the blowups of $M$ along the three sheaves in the equation.
\m\m \centerline {4. Properness of $U\to Y_0$}
\m\m 
If it is not known whether  $U\to Y$ is proper, instead
let $f:Y_0\to M$ be the universal singular manifold over $M$  such that
at every point $p\in Y_0$ there is an index $r$ such that 
$f^*{\cal F}_i/torsion$ is principal in a neighbourhood of $p$ when
$r\vert i.$

\m {\bf 6. Corollary.}  The map $U\to Y_0$ is always proper. Therefore
$U\to M$ is proper if and only if $Y_0\to M$ is proper.
\m Proof. The equation above shows
that the local isomorphism type of the pullback modulo torsion
of ${\cal F}_{(n+2)^m-1}$ stabilizes as $m$ increases. This does immediately
imply the corollary without further work, but a less mysterious proof
 is by considering arc lifting
as we shall do in section 8.
\vfill\eject

\centerline{5. The description of $U$ as a graph}
\m
\m {\bf 7. Theorem.} Suppose $M$ is projective with very ample ${\cal F}$  and that $Y_0\to M$ is proper. 
Let $W=Proj\oplus_i\Gamma({\cal F}_i).$ Then $W$ is an algebraic
variety and   $Y_0\to M$ is not only a proper map, it is actualy
 the universal solution of
resolving the indeterminacies of the rational map $M - \to W,$ ie it 
is the closure of the `graph' of the rational map $M- \to W.$ 
\m The following is an immediate consequence of theorems 5. and 7. 
\m  {\bf 8. Corollary.} The map $M- \to W$ is a morphism if and only if 
$M$ is nonsingular. 
\m {\bf 9. Remark.} Assume that $M$ is a singular projective variety. Then
 the singular manifold $W$ consists of a single point
if and only if $M$ is a linear projective space and $H$ is a hyperplane.
\m This shows that $dim(W)$ need not {\it always}  be as large as $n.$ 
It is of course bounded above by the Iitaka dimension of $K+(n+1)H.$ 
We have to be careful if we want to work
 by induction on Iitaka dimenson
as
for any Gorenstein variety besides projective space
which has
isolated irrational singularities,  
$K+(n+1)H$ is very ample (see [16] theorem 8.8.5).
Beginning in section 7 we'll consider the consequences
of forcing the calculation when $n$ is intentionally chosen
to be different than the dimension of $M.$ 

\vfill\eject \centerline{6. The terms of degree ${{(n+2)^i-1}\over{n+1}} .$ }
\m\m The terms ${\cal F}_i$ where $i$ is a partial sum of a geometric
series to the base of $n+2$  play a particular role, as if I replace ${\cal F}$ by 
such an ${\cal F}_i$ then the auxiliary singular manifold $Y_0$ is unaffected,
and the effect on the sequence ${\cal F}_0, {\cal F}_1,...$ is
to `truncate' by passing to the subsequence of multiples
of a power $(n+2)^j.$ 
\m  
Properness of $Y_0\to M$ is equivalent $\oplus {\cal F}_i,$ 
being (=locally) of finite type. That is to say,
on each bounded open subset of $M$  
 ${\cal F}_{(n+2)^i}^{\otimes (n+2)}
\to {\cal F}_{(n+2)^{i+2}}$ is a surjective map of rank one
coherent sheaves for large $i.$ Once this happens for one value of
$i$ then it happens for all sufficiently large values of $i.$ 
Surjectivity of the map from the $n+2$
symmetric power of the terms of degree $1$ to
terms of degree $n+2$ is neither necessary nor sufficient for  nonsingularity
of $Bl_{\cal F}(M),$
but it follows nevertheless therefore from 
 Hironaka's theorem [3] and [20] Theorem 3.45 (see also crucial
references therein)  that on a bounded open set once ${\cal F}$ is replaced by such an ${\cal F}_i$ that
 there is  always a choice of ${\cal F}$ such that
$\oplus_i {\cal F}_i$ has generators of  degree one. Conversely when ${\cal F}_1$
generates, then ${\cal F}\otimes{\cal F}_1/torsion$ is a resolving sheaf.

\vfill\eject
\centerline{7. Filtration of the homogeneous coordinate ring of $V$ }
\m\m The results above
provide a functorial relative affinization of Nash's manifold 
which is analogous to the intrinsic relative affinization of $V.$
This is not precisely automatic, for example no functorial affinization
of Hironaka's and 
Spivakovsky's model [14] is known. 
\m
If $n$ is taken to be $1+dim(M)$ in the definition of
the multiplication (but not in the exterior degrees) then a 
series of subsheaves of the ${\cal F}_i$ arises.  It is confusing to say `let us no longer
assume that $n$ is equal to the dimension of $M$' and so
let's instead continue to let $n=dim(M)$ and use a contrivance: for each
integer $N\ge -1,$ whenever  $j$
is a power $j=(N+2)^i$ let
$$X_j[N]={\cal F}^{\otimes(n+1){{(N+2)^i-(n+2)^i}\over{N-n}}}\otimes {\cal F}_{ (n+2)^{i-1}}\otimes{\cal F}_{(n+2)^{i-2}}^{\otimes N+2}
\otimes {\cal F}_{(n+2)^{i-3}}^{\otimes (N+2)^2}\otimes...
\otimes{\cal F}_1^{\otimes (N+2)^{i-1}}/torsion$$ depending on $N,$ $n$ and $j.$
The exponent of ${\cal F}$ is positive regardless of the relative magnitude
of $N$ and $n.$  If $j$ is not a power of $(N+2)$ 
then define $X_j$ to be a product according to the base $N+2$ expansion of $j,$
and let
$${\cal F}_j[N]=lim_s {\cal H}om(X_j[N]^{n+s},X_j[N]^{N+s}\otimes{\cal F}_k/torsion)$$ where $k$ results by replacing $N$ by $n$ in the base $N+2$ 
expansion of $j.$ (When $N=-1$ and base $N+2$ expansions are ambiguous 
also take the colimit over expansions).
Taking $N=n$ we have 
$${\cal F}_j \subset {\cal F}_j[n]$$
an integral map. 
For $N$ taking values besides $n$ then
the ${\cal F}_j$ are nested
$${\cal F}_j[n+1]\subset {\cal F}_j[n]\subset ...\subset{\cal F}_j[-1]$$
and it is also easy to see that for each fixed $n$ the $\oplus_j {\cal F}_j[N]$
form a graded ring sheaf.   We will also  abbreviate $\oplus_i{\cal F}_i[N]$
by the symbol ${\cal F}[N].$ 
\m {\bf 10. Theorem.} 
As
$N$ passes along the sequence $-1,0,1,2,...,n+1$ the ${\cal F}[N]$
describe then a filtration of the sheaf of  homogeneous coordinate rings of $V.$
\m Proof. By 14 d) ${\cal F}_j[-1]\subset \tau_*{\cal O}_U(jK_U).$
\vfill\eject\noindent  For each $N$
let
$f:Y_0[N]\to M$ be the universal singular manifold such that each point $p$
has an index $r$ such  that $f^*{\cal F}_i[N]$ is principal in a neighbourhood of $p$ for all $r | i.$  The resulting $Y_0[N]$ is independent of the choice of ${\cal F}$ if ${\cal F}$ is
invertible, and let us  take ${\cal F}={\cal O}_M.$ It will follow
from Theorem 14
 
\m
 {\bf 11. Theorem.}
\begin{enumerate}
\item [a)] When $N=-1,$ the integral closure of  
${\cal F}[N]$ is the homogeneous
coordinate ring of $V,$ and so the  $Y_0[-1] -\to V$ is the inverse
of a finite (=proper locally affine) morphism. 
\item [b)] When $N=n,$   $Y_0[N]$ admits a finite  map
 to $Y_0$ over $M.$
\item [c)] When $N\ge n+1,$ $Y_0[N]$ is isomorphic to $U$.
\item [d)] For values of $N$ between $-1 $ and $n+1$ there are
meromorphic maps whose inverses are carried by morphisms of the
affinizations. 
\item [e)] The first  map $Y_0[n+1]- \to Y_0[n]$ is proper and  holomorphic.
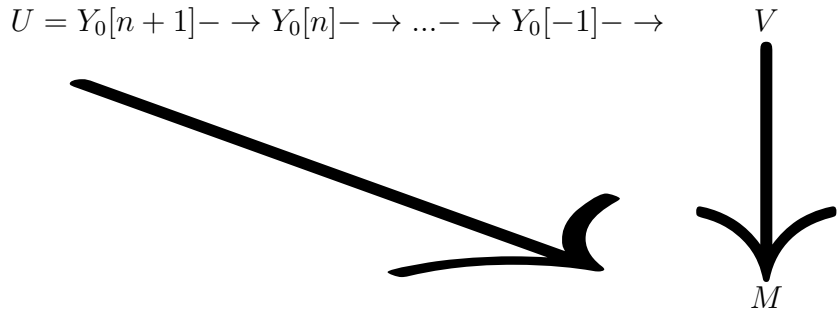
\m \m $$\begin{matrix}U=Y_0[n+1] - \to Y_0[n] - \to ... - \to Y_0[-1]  -\to  &V \cr
\resizebox{8cm}{2cm}{$\searrow$}&\resizebox{2cm}{2.5cm}{$\downarrow$}\cr
& M\end{matrix}$$
\end{enumerate}
\m
It holds from [13,21,22,23] that the map $V\to M$ is 
a proper morphism if $M$ is locally algebraic (this is full
generality if the singularities of $M$ are isolated [5]).  
The diagonal vertical arrow and the vertical arrow
are the two models $U$ and $V.$ The existence of a chain of proper
maps in either direction between $U$ and $V$ over $M$ is not a new
open question. Unless both were false, it is equivalent
to the logical equivalence between the questions considered
separately of the  properness of $U$ over $M$
and the properness of $V$ over $M.$ 
\vfill\eject\centerline{8. Nash Arcs}
\m Pulling back $\oplus_j{\cal F}_j$ modulo torsion along  an analytic arc
$t\mapsto \gamma(t)\in M$ which is defined on a domain of times
$T\subset {\mathbb C}$ 
gives a graded ring of  coherent sheaves 
 on a domain in the complex number line. Suppose that the arc contains
a regular point. Then the sheaf is finitely generated except at
a discrete set of times $t,$  so there is no harm assuming that
$\gamma(t)$ is a regular (smooth) point except when $t=0.$ 
The restriction of  $\gamma^*{\cal F}_{(n+2)^i}^{-n-2}\cdot
\gamma^*{\cal F}_{(n+2)^{i+1}}$ to the open set $T\setminus \{0\}$
is canonically isomorphic to ${\cal O}_{T\setminus \{0\}}$ and contains
a canonical copy to ${\cal O}_T.$ 
On $T$ there is a thereby distinguished natural isomorphism of
the torsion free part to
 $t^{-d_i}{\cal O}_{T}$ for a natural number $d_i$ depending
on $i.$ Thus we obtain a sequence of 
natural numbers $d_1, d_2,...$ 
\m {\bf 12. Theorem.} The sum $\sum_{i=1}^\infty d_i$ is finite
if and only if $\gamma$ lifts to $U$, and it is then  equal to 
the local intersection product $(K_U-K)\cdot \gamma$
where
 $K$ is a Cartier divisor such that ${\cal O}_U(K)$ is the pullback
mod torsion of $\Lambda^n\Omega_M.$ In other words, passing to
the local analytic ring ${\mathbb C}\{\{t\}\}$ at the origin,  
 the minimum number
of generators $gen(\gamma)$ of the local 
algebra     $\gamma^*(\oplus_i {\cal F}_i)_{\{0\}}$ 
   over
 ${\mathbb C}\{\{t\}\}$  
is  given by the equation
$$\boxed{gen(\gamma)=(K_U -  K)\cdot \gamma + 1}$$
\m Proof. The integral
extension from section 10
  $\oplus_i{\cal F}_i\subset \oplus_i{\cal O}(D_i)$ 
pulls back modulo torsion to an isomorphism.  
Then $d_0=1$ and  
$d_i=\gamma\cdot (D_{(n+2)^i}-(n+2)D_{(n+2)^{i-1}})$
for $i\ge 1.$ By the definition of the $D_j$ this equals $\gamma\cdot (K_i-K_{i-1}).$   
Adding over $i=1,2,...$ yields a
 collapsing sum adding to $\gamma\cdot(K_U-K).$ 
\m For example if $C$ crosses $K_U-K$ transversely at one point
then $\gamma^*\oplus{\cal F}_i$   is isomorphic to a graded sheaf
of algebras generated by $t$ in degree  $1$ and
 $t^{(n+2)^s-1}$  in degree $(n+2)^s$
where $s$ is, I believe,
the number of Nash blowups needed to make the intersection transverse.
The intersection number is then the local winding number of $T\setminus \{0\}$
about one of the components of $K_U-K.$

\m 
\vfill\eject\centerline {9. Examples and discussion.}
\m\m This section will have no theorems, but is included for context.
In the case when ${\cal F}$ is very ample the ${\cal F}[N]$ are generated
by global sections as long as $N\ge n.$ 
The square brackets instead of round brackets around $N$ are to
avoid confusion with the notation
of level
 in modular forms. For
example  if  $M$ is a modular
Riemann surface of level $G \subset {\mathbb P}Sl_2({\mathbb Z}) $ a torsion free group,
and ${\cal F}$ is the section sheaf of a line bundle, including one 
global section which crosses  each cusp transversely,
 then in the inclusions of four graded rings
$$\Gamma  
{\cal F}[2]\subset \Gamma{\cal F}[1]
\subset \Gamma {\cal F}[0]
\subset \Gamma  {\cal F}[-1]$$
all four rings are the same; they are all equal  to the ring of modular forms of even weight and
level $G.$ The composite of the four equalities induces the composite of the Nash resolution 
with the inverse of the relative canonical morphism.
An example where $N$ matters is the  case
of something higher-dimensional such as the
example of Kollar and Ishii [18],
the solution set $M$ of the equation
$x_1^3+...+x_4^3+x_5^6=0,$ 
and let us here repeat the beginning part of their discussion. Here $V\to M$ is
an isomorphism because $M$ already has canonical singularities. These
can be resolved by blowing up reduced points. First blowing up the singular point
yields as exceptional component the projectivized tangent cone at the
origin, a cone on a cubic hypersurface in ${\mathbb P}^3,$ 
 which is ruled. Nevertheless it must occur
in any resolution because of Nash's argument [15] about 
irreducible components in the space of arcs.
The  blow-up in turn of the cone point 
of the projective subvariety  provides a resolution $V'\to V$ with
  a second component, a full cubic hypersurface in ${\mathbb P}^4$
 which is not an arc  component because
first infinitesimal neighbourhoods of lines split;
but is now essential because it is non-generically ruled.  Thus
the paper produces two essential components which live `above'
the smallest canonical singularities model $V$ but for different
reasons. 
Also [18] says that even if we did not now that 
the first component were an arc component, and even though 
it is not one of the essential crepant components of [7], its
essentialness also follows
 because of having minimal discrepancy.
\vfill\eject\noindent 
In this case 
  $\Lambda^4\Omega_V/torsion ={1\over{x_5^5}}
I
dx_1dx_2dx_3dx_4$ has invertible reflexivication
$ {1\over{x_5^5}}{\cal O}_Vdx_1dx_2dx_3dx_4$  
where $I=(x_1^2,x_2^2,...,x_5^5).$ The section  ${1\over{x_5^5}}dx_1...dx_4$ 
of the canonical line bundle of the resolution a  has a simple
zero at the first exceptional component and a double zero at the second
one.
Taking ${\cal F}$ to be the unit ideal, each  $ {\cal F}_i
$ contains $I^i {1\over{x_5^{5i}}}(dx_1...dx_4)^i$ and is of the
form $I_i{1\over{x_5^{5i}}}(dx_1...dx_4)^i$ for $I_i$ a suitable
adjunction ideal.
Defining the sequence $g_1,...,g_5=x_1^2,...,x_5^5,$ 
then ${\cal F}_1=I{1\over{x_5^5}}dx_1...dx_4$ is ordinary $4-$ forms, and ${\cal F}_6$ 
contains ${\cal F}_1^6$ but is slightly larger; it
is spanned over ${\cal O}_V$ by the alternating forms
$$g_0dg_1dg_2dg_3dg_4-g_1dg_0dg_2dg_3dg_4+g_2dg_0dg_1dg_3dg_4$$ $$-g_3dg_0dg_1dg_2dg_4+
g_4dg_0dg_1dg_2dg_3$$
for $g\in {\cal F}{\cal F}_1={\cal F}_1.$ In other words
symmetric linear six-forms in alternating differential four-forms 
are 
 a special case of an alternating sum of products of
differential four-forms against four-forms in four-forms.
Amusingly, the expression
above can be rewritten as $g_0^5 d(g_1/g_0)...d(g_4/g_0)$ showing that 
this expression
is then antisymmetric in the $g_i.$ This same expression 
which describes a typical $n$-form on a coordinate chart with
an extra zero of degree $n+1$ on the exceptional locus 
describes also a typical element in the image of 
the universal connection.  
Let us  explain this
and also something more general. 
\m {\bf 13. Observation} The expression $g_0^{n+1}d(g_1/g_0)...d(g_n/g_0)$ is
antisymmetric under permutations of $g_0,...,g_n$ and results
from applying   
the universal ${\mathbb C}$-linear connection $\nabla$
to $dg_0...dg_n$ using Leibniz rule. Therefore ${\cal F}_6$ is generated
by the image of $\nabla$ if $g_0,...,g_4$ run over sections of the
product sheaf ${\cal F}{\cal F}_1.$   
\m Just generally, if ${\cal G}$ is  torsion-free coherent of rank
$s$ then  
$$
0 \ \buildrel H \over \rightarrow  
\ \Lambda^{s(n+1)}{\cal P}{\cal G}\otimes S^j({\cal G})/torsion 
\ \buildrel H \over \rightarrow 
\ \Lambda^{s(n+1)-1}{\cal P}{\cal G}\otimes S^{j+1}({\cal G})/torsion
  \buildrel H \over \rightarrow ... $$
an  exact complex for $j\in {\mathbb N}.$ 
\vfill\eject\noindent
The sheaf ${\cal P}({\cal G})$ is the kernel of $q:\Omega_{V({\cal G})}\to \Omega_M,$
(that is, one-forms which are zero on $M$)  pulled back as a coherent sheaf along
the inclusion $i:M\subset V({\cal G})$ of the zero section, to arrive at  $i^*Kernel(q)/torsion.$ If $y$
is a section of ${\cal G}$ viewed as a linear function on the fibers
of $f,$ it is already in $Kernel(q),$  and for $x$  a  local holomorphic function on $M,$ the connection $\nabla$
is merely the deRham extending the deRham differential $d$ on $M$ and satisfying
$\nabla(xy)=x\nabla(y)+ydx.$  
\m
Now pulling back ${\cal P}({\cal G})$ further along $f:V({\cal F})\to M$
gives $f^*{\cal P}({\cal G}),$ and $H$ and $\nabla$ extend by Leibniz rule.
$H$ is a 
contracting homotopy; the formulas
$$H\circ \nabla + \nabla \circ H = (s(n+1)-j)$$
$$ H \circ H = 0$$
$$\nabla \circ \nabla=0.$$
$H$  commutes with $x$ and sends $dy\mapsto y \mapsto 0.$
\m
The result $f^*{\cal P}({\cal G})/torsion $ is the 
direct sum of the sequences shown above 
plus  a finite
number of incomplete parts of sequences
exact except possibly
at the leftmost term, which is 
$\Lambda^i({\cal G}\otimes_{{\cal O}_M} \Omega_M) $ when $i$ is
small enough that this is nonzero.
The Atiyah sequence $H:\Lambda^1{\cal P}({\cal G})\to S^1 {\cal G} \to 0$ 
with its ${\mathbb C}$-linear
splitting $\nabla:{\cal G}\to {\cal P}(G)$ is one of those incomplete
parts, the one for $i=1,$  although the information there determines the other sequences.
The Atiyah sequence itself does not transform well since the kernel
term $\Omega\otimes {\cal G}/torsion $ is
nonvanishing on the zero section $M.$ 
The notion of Atiyah is that although this sequence is exact, still
$K$-theoretic information can be extracted as Whitehead torsion.  
The same should be true of the exact sequences which transform better too.
That is, the exact sequences of sheaves whose sections vanish on the
zero section $M$ are not affected by various modification of 
the zero section such
as intermediate blowups between $M$ and $Bl_{\cal G}(M)$.
Though algebraic cycles representing the Whitehead torsion might
change dimension reminiscent of dimension changes in arc components.

\vfill\eject\noindent  Esnault, Viehweg and 
Verdier show in [12], Appendix B,  how to extract higher Chern classes using
Deligne's concept of connections with logarithmic poles; and they
mention that the assumption of a 
normal crossing divisor of multiplicity one is not
needed. In place of the composite $\Gamma_i^{\alpha_i}$ there, 
one could try looking for algebraic cycles in the higher
exterior power complexes here rather than only iterating
the morphism.  Here, where ${\cal G}$ has rank one,
this is just the fact that $\nabla \circ H$ is multiplication by $n+1$ on
$\Lambda^{n+1}{\cal P}{\cal G}.$ It means that $H$ is an 
embedding of $\Lambda^{n+1}{\cal P}{\cal G}/torsion$ in 
$\Lambda^n{\cal P}({\cal G})\otimes {\cal G}/torsion,$ and in fact the
exact sequence going further to the right (with differential $H$) is exactly like the
sequence used in the theory of Castelnuovo-Mumford regularity;
it is a torsion-free and exact Koszul complex twisted by 
tensor powers of ${\cal G},$ and
 it is exact even though ${\cal G}$ is an arbitrary
torsion free coherent sheaf of rank one. Generically the image
is equal to the rank one subsheaf $\Lambda^n(\Omega{\cal G}) \otimes {\cal G}/torsion=(\Lambda^n\Omega)\otimes {\cal G}^{n+1}/torsion;$ in fact
it is a bit larger. The inclusion composed with the inverse of the
isomorphism of $H$ onto its image is
precisely the map of theorem 1  
for the case ${\cal F}={\cal O}_M.$

\m 
We see from the observation that if we coordinatize the blowup of $(g_0,...,g_{n+1})$ 
we are adjoining an $n$ form with extra poles of degree $n+1$ on the exceptional
locus. Without blowing up the ideal, the presence of the additional
generators means 
that we are enlarging  the adjunction ideal by an inclusion $I^{6}\subset I_6$
without affecting the codimension one primary components (there are none);
elements of ${\cal F}_6$ are not six-fold symmetric powers of one-forms anymore,
they can be viewed as having additional poles on $V$ at primary ideals whose associated prime
is the maximal ideal of the cone point. 
\m If $\pi:V'\to V$ is a resolution the sheaf $\oplus_j \pi_*(\omega_{V'}^{\otimes j})$
is just $\oplus_j \omega_V^{\otimes j},$ and any isomorphism $\omega_V\to {\cal O}_V$
gives an isomomrphism between this and (the sheafification of) the polynomial
algebra ${\cal O}_V[T].$ The $ {\cal F}[N]$  describe a filtration of the
graded algebra and we have noticed that ${\cal F}_6[4]$ is like six-fold symmetric
powers of alternating differential 4-forms with extra poles of codimension four,
 it is just an integral extension of the un-adorned ${\cal F}_6$;  
while there must exist a number $r$ such that the  
integral closure of ${\cal F}_j[-1]$ is  all of 
$ ({\cal F}_1^{\otimes j})^{**} =\omega_V^{\otimes j}$ when $r|j.$

\m \vfill\eject\centerline{10. Divisorial approximation}
\m\m Let us say that an inclusion ${\cal F}\to {\cal G}$
of torsion-free coherent sheaves of rank one
is {\it integral} if for each  Stein open set $U$ there is some ideal sheaf 
${\cal Y}$ so that
$\Gamma(U,{\cal YF})  = \Gamma(U, {\cal YG}).$  

\m The construction of the relative canonical model by taking direct images suggests
an analagous construction, and up to integral morphisms this is possible.
\m Let $K_i$ be the Cartier divisor on  $U$  
which is the pullback modulo torsion of $\Lambda^n\Omega_{M_i}.$
Fix a natural number $N\ge -1.$ 
Let $\tau:U\to M$ be the structural map. Define  Cartier divisors $D_1, D_2,...$
on $U$ to be the integer linear combinations of the $K_i$ given  by the equation
$$D_{(N+2)^j}= K_j+(N+1)(K_{j-1}+(N+2)K_{j-2}
+(N+2)^2K_{j-3}+...+(N+2)^{j-1}K_0)$$ and when $i=a_0+a_1(N+2)+...+a_s(N+2)^s$
with $0\le a_i<(N+2)$ let $$D_i=a_0D_1+a_1D_{(N+2)}+...+a_sD_{(N+2)^s.}$$ 
\m If there is any ambiguity about the choice of $N$ we will
write $D_i[N]$ to denote the relevant number $N$ which was used
in constructing $D.$
\vfill\eject\noindent {\bf 14. Theorem.}  Let ${\cal F}={\cal O}_M.$
\begin{enumerate}\item [a)]  For each  $N$ and $j$   There is a natural
map of torsion-free coherent sheaves of rank one
$${\cal F}_{j}(N)\to \tau_*{\cal O}_U (D_j[N]).$$
It is an integral map for all $N$ and $j.$ 
Therefore  $ {\cal F}[N] \subset
\oplus_i \tau_*{\cal O}_U(D_i[N])$ is an integral map of ${\cal O}_M$-algebra
sheaves and, more strongly, of underlying coherent sheaves.
\item [b)] For each $N$ and each sequence of numbers increasing
with divisibility $i_1 \vert i_2 \vert i_3| ...$
the sequence of ${\mathbb Q}-$ divisors on $U$ is an ascending
sequence
$${1\over{i_1}}D_{i_1}[N]\le  {1\over {i_2}}D_{i_2}[N]\le {1\over {i_3}}D_{i_3}[N], ....$$ 
and on any bounded subset of $Y_0[N]$ it is a finite ascending series
of ${\mathbb Q}$-divisors. Therefore for each $N$ there is a limiting ${\mathbb Q}$-divisor
$D[N]$ on $U$ such that at each point of $U$ there
is an index $r$ such that $D_i[N]=i\cdot D[N] $ whenever $r|i.$
\item [c)] For $N=-1$ the divisor $D[N]$ is merely 
the canonical divisor $K_U$ of $U.$ That is, $D[-1]=K_U.$
\item [d)] Let $\pi:M'\to M$ be any resolution.
For each $i$ the pushforward $\tau_*{\cal O}_U(iK_U)$
is the same as $\pi_*{\cal O}_M(iK_M).$ 
We can therefore use $U$ in place of $M'$  and $D[-1]$ in place
of $K_{M'} $  in defining the
homogeneous coordinate ring sheaf  of $V.$
\item [e)] Recall $D_i=D_i[n].$   
The fibers of $U\to Y_0$ are (complete) projective varieties.
Let $C$ be a complete curve in $U.$  Then $C$ is contained in just one fiber of $U\to Y_0$ if and only if
the right side of the equation 
$$\boxed{K_U\cdot C= -(n+1)\cdot \sum_{i=0}^\infty D_{(n+2)^i}\cdot C}$$ 
has only finitely many nonzero terms.  When $C$ is not contained
in a fiber the right side is eventually a geometric series. In
every  case 
the equation remains valid formally; i.e. as the sum of 
a convergent series of 
$p$-adic numbers for any prime divisor $p$  of $n+2$. 
\end{enumerate}
\vfill\eject\noindent 
{\bf 15. Corollary.} The `discrepancy' $K_U\cdot C - K\cdot C$ 
is the sum over points $p\in C$ of the $gen_p(C)-1$ where $gen_p$ is the local generating number
at $p.$ This differs from the  actual discrepancy by $(\tau^*K_M-K)\cdot C.$ I think that the sum is the dimension of   the vector space
$(i^*{\cal F}/torsion) \otimes_ {\cal F}{\cal O}_M$
if $i:C\to U$ is the map. If $\gamma$ does not lift, the cardinal
dimension of the vector space does not carry full information; the
Grothendieck class of the module must be represented by 
a Poincare series or a $p$-adic number in that case.

\m 
Proof of d). Assume $M$ is contained in a disk. 
Choose a global section of the $i$'th 
power of the  canonical sheaf of the resolution $M'.$
Choose a point of $U,$ and choose
a stage of the Nash tower $M_m\to M$ which contains
a neighbourhood $T$ of this point. The form on the regular
locus of $M'$ pulls back to a form on a resolution of the pullback of
$M'$ and $M_m$ over $M$ and forward by an isomorphism between an open
set and the complement of a codimension two subset of $T,$ and then 
extends across all of $T.$  For the converse
use Hilbert's basis theorem for subsheaves of 
of $(\Lambda^n\Omega_M^{\otimes i})^{**}.$ Replacing $M$ by a bounded
open subset if necessary, there is 
 a number $m$ so   $\tau_*{\cal O}_U(iK_U)
=\tau_*{\cal O}_U(iK_m).$ 
In the case $M_{m+1}$
is normal, since ${\cal O}_U(iK_m)$  agrees on the regular locus of $M_{m+1}$ with  the pullback modulo torsion of $\Lambda^n\Omega_{M_m}^{\otimes i},$ we have a section of
a locally free sheaf,  which then extends across the codimension-at-least-two
singular locus of $M_{m+1}.$ 
The extended section belongs to the larger sheaf $\Lambda^n\Omega_{M_{m+1}}^{\otimes i}/torsion$ 
and  does  pull back 
to a section of $\omega^i$ on a resolution.
If $M_{m+1}$ is not normal, boundedness of 
 $M_{m+1}$ in the larger  manifold implies that
there is an upper bound on the number of initial consecutive finite and nontrivial Nash blowups of any open subset of $M_{m+1}.$ A higher
Nash blowup $M_N\to M_{m+1}$ yields, by deleting an exceptional set $E\subset M_N$  and its image $C\subset M_{m+1}$ of codimension
at least  two, a finite map which is a resolution.
Our section restricts on $M_N\setminus E \subset U$ to a section of the pullback
of $\Lambda^n\Omega_{M_m}^{\otimes i}$ modulo torsion. 
Except on the inverse image of $C$ 
this agrees  with a section of the
 of the (locally free) pullback modulo torsion
of $\Lambda^n\Omega_{M_m}^{\otimes i}$ to the normalization of $M_{m+1}.$ 
It extends to a global section
as in the normal case.
\vfill\eject\noindent
It is tempting to summarize the results by saying that 
$D[n+1]$ is relatively ample on $U,$ $D[n]$ is relatively
eventually basepoint free, and $D[-1]$  is canonical with classical 
relative plurigenera (the same as a resolution). 
A difficulty is that we have to clarify what `ample' and
`eventual' actually mean, because both definitions implicitly
refer to an index  $r(p)$ associated to each point $p\in U$.
In the algebraic case semicontinuity of the best $r(p)$  in the Zariski
topology and quasicompactness give a bound and it had not been
necessary to distinguish the cases. Here
we need to clarify whether we intend the $r(p)$ to be
bounded above.
The assumption that $M$ is a manifold is not important since
any singular manifold is contained in a nonsingular manifold. 
\m {\bf 16. Definition (clarification).} 
Let $\tau :U\to M$ be a map of smooth complex manifolds.
\begin{enumerate}
\item [a)] A Cartier divisor (or ${\mathbb Q}$-divisor)
$D$ on  $U$ will be called  {\it uniformly  relatively ample} 
if some  multiple $rD$ is relatively very ample. 
\item [b)] It will be called  {\it non-uniformly relatively  ample} if
 nevertheless for  each point
$p\in U$ there is an  index $r(p)$ depending on $p$ 
and neighbourhoods $T$ of $p$ and $S$ of $\tau(p)$
such that 
 $\Gamma(S, \tau_*{\cal O}_U(rD))$ 
separates points and tangent vectors on 
$T;$
 ie that  
the blowup  $Bl_{\tau_*{\cal O}_U(r(p)\cdot D)}\to M$ 
contains a copy of $T$ and the restriction of $\tau$ to $T$
is induced by the structural map of the blowup. 
\item [c)] It will be called {\it uniformly eventually relatively
basepoint free} if there is an index $r$ such that ${\cal O}_U(rD)$
is relatively basepoint free.
\item [d)] It will be called {\it non-uniformly eventually relatively
basepoint free} if nevertheless every point $p$ has in index $r(p)$ depending
on $p$ and  neighbourhoods $T$ of $p$  and $S$ of $\tau(p)$ such
that the restriction map 
$\Gamma(\tau^{-1}(S), {\cal O}_U(rD))\otimes {\cal O}_T \to {\cal O}_T(rD))$
is surjective.\end{enumerate}
\m The terminology is chosen so that the uniform case is included in 
the uniform case, characterised by the boundedness of the $r(p).$ 
\vfill\eject\noindent
In this sense then 
\m {\bf 17. Summary.} The  
Cartier ${\mathbb Q}$-divisors $D[N]$ on the manifold $U$ are such
that 
$$D[N] \hbox{ is } 
\left\{\begin{matrix}
\hbox{ non-uniformly relatively ample,}& N\ge n+1\cr
\hbox{ non-uniformly relatively eventually basepoint free}, & N\ge n\cr
\hbox {canonical with classical relative plurigenera},& N=-1
\end{matrix}\right.$$

\m Things can be made absolute rather than relative
in the case $M$ is quasi-projective. We let ${\cal F}={\cal O}_V(H)$
for $H$ a hyperplane section. Then,  for example
\m {\bf  18. Theorem} (some global variants) Suppose $M$ is quasi-projective and $H$ is
very ample on $M.$ 
  Then  
\begin{enumerate} \item [a)] Each divisor $D_i+(n+2)^i(n+1)\tau^*H$ is
absolutely basepoint free on $U.$
\item  [b)]  For each $N\ge n+1$  
and every point $p\in U$ a suitable integer multiple
of the ${\mathbb Q}$-divisor $D[N]+(n+1)\tau^*H$ 
defines a projective embedding of a neighbourhood of $p.$
\item  [c)]
For $N=n$ a suitable integer multiple of $D[N]+(n+1)\tau^*H$ depending on $p$ 
has base locus disjoint from $p$
\item  [d)] For $N=-1$ and for all $j$ the
global sections on $U$ $\Gamma({\cal O}_U(j(D[-1]+(n+1)\tau^*H))$ are the just
 the global sections  $\Gamma(M',j(K_{M'}+(n+1)\pi^*H)$ for
a resolution $\pi:M'\to M.$
\item [e)] If $M$ is normal, the fiber in $U$   over each singular point of $M$ always contains at least
one contractible complete curve $C$ with $D[n]\cdot C=0.$ 
\m Part d follows  from Theorem 14  by the projection formula.

\end{enumerate}

\m {\bf 19. Remark.} It appears also  to be  true in the toric case
if one makes the strange assumption that all
the  $\Lambda^n\Omega_{M_i}/torsion $ on the terms of the Nash
tower 
are  reflexive then for some reason the basepoint-free results
extend all the way down to $D[0].$ 

\vfill\eject \centerline {11. Questions}
\m\m  The obvious idea, before trying calculations, is to look for
successively stronger basepoint free theorems to attempt to find
a series of proper maps in either direction relating $Y_0$ with $V,$ 
following Mori's discovery [13]. If $M$ locally admits a smooth one-dimensional
foliation, for example, then the next map $Y_0[n]\to Y_0[n-1]$ is
also a locally projective holomorphic map. Note that the locally
projective morphisms are in the opposite direction from the affine
morphisms. In general there is no 
idea what modification of $Y_0[n]$ 
can creates a useful relatively-basepoint-free divisor.   
Theorem 1 is
merely tautological, as opposed to more substantial basepoint free
theorems in algebraic geometry. 

\m If $M$ is not only a complex manifold but a singularly foliated
complex manifold then $n$ may be replaced not by the codimension but
the dimension of the foliation, and all the above holds, except
it refers not to the Nash manifold, but to the smoothly foliated
singular manifold which results by repeatedly blowing up the foliation
itself.
\m {\bf 20. Conjecture.} Suppose $M$ is a complex vector space
singularly foliated by the faithful split action of a commutative
Lie algebra. Then the foliation is resolvable by a locally projective
degree one morphism from a singular variety with a nonsingular foliation
if and only if the set of roots (ignoring multiplicity) forms a basis
of the dual of the Lie algebra.
\m This is worked out in the case of toric resolutions in an unpublished arXiv
preprint [19]. The same preprint claims without proof
\m {\bf 21. Theorem.} Suppose $M$ is a complex singular projective variety
with very ample divsor $H$ and that $M$ has a singular foliation of dimension $s$
such that the canonical divisor $K$ of the foliation satisfies
that $K+(s+1)H$ is a finitely-generated divisor of Iitaka dimension one. Then the foliation can
be resolved by one or more Nash blowups of the foliation.
\vfill\eject\noindent This is  true because the successively higher Nash blowups
are induced by maps to algebraic curves $W_i.$  
The rational function fields ${\mathbb C}(W_i)$ 
are all contained in a subfield of a field of transcendence degree one.
Therefore the rational function fields of the $W_i$ stabilize and then
the $W_i$ are bounded by the unique birational model. 
\m Note also that  even though the
ring $\oplus_i \Gamma({\cal F}_i)$  may  conceivably fail to be finite type,
the sheaf $\oplus_i{\cal F}_i$ actually must then be finite type. That is, 
since the  $Bl_{{\cal F}_i} M$ are eventually isomorphic over $M$ for
large $i,$  the ${\cal F}_i$
pull back to invertible sheaves  on suitably high Nash blowups.
By Theorem 1 then when $i$ is a large
power of $N+2,$ ${\cal F}_i^{N+2}\to {\cal F}_{i(N+2)}$ is onto.

\m The canonical sheaf is the simplest
primary component of the highest exterior power of the differentials.
This component is invertible precisely when $-K_0$ is an effective exceptional
divisor in the Nash tower.
In general, if it is not invertible, if we blow up only this  primary component
the resulting map is not the identity; if we do so repeatedly
\m {\bf 22. Conjecture.} 
If we blow-up only upon the
codimension one primary component of $\Lambda^n\Omega/torsion$  
at each stage,ie its reflexivication,  then after finitely many steps 
the reflexivication will be invertible; there results
a proper morphism from a singular manifold with invertible canonical sheaf.
\m The conjecture is easily verified to be true for toric surfaces.
The pullback modulo torsion of the canonical sheaf
 is an invertible  subsheaf
of the reflexive canonical sheaf,  
and so there is an effective 
Weil ramification divisor.
 The process finishes if and only
if the Weil ramification becomes zero and that is what the conjecture
asserts.
\vfill\eject\m Hironaka's and Spivakovsky's 
theorem is that for surfaces  normalized Nash blowups do yield
a proper map of a smooth manifold to $M$. There appears to be
a  great distance between the normalized and non-normalized Nash blowups.
\m When ${\cal F}$ is taken to be ${\cal O}_M$ then 
functoriality in the title is with repsect to degree-one  holomorphic
maps of $M.$ Strictly speaking it applies only to ${\cal F}[N]$ for
$N\ge n,$ after passing to integral closure. For example, the Grauert-Riemannschneider sheaf is the integral closure of ${\cal F}_1[-1]$ when ${\cal F}={\cal O}_M$, and this is not a functor. The failure of functoriality however is
only due to the fact that ${\cal H}om$ is contravariant in one of its
arguments.
\m That is, a corollary of part d) of Theorem 14 and the definition
of $X_i[-1],$  yields an elementary calculation of 
Grauert-Riemenschneider's sheaf whenever ${\cal F}$ is invertible
(e.g. if ${\cal F}={\cal O}_M$)
\m {\bf 23. Theorem.} Suppose $M$ is bounded (within a possibly larger manifold
of the same dimension). 
The Grauert-Riemenschneider sheaf of $M$
is the integral closure of 
   ${\cal F}^{-n-1}{\cal H}om( {\cal F}_{(n+2)^i-1},
{\cal F}_{(n+2)^i})$
for all sufficently large value of $i.$
\m The calculation of the sheaf has already been  within the range
of computer since resolution of singularities algorithms exist.
A simpler method would be to use the formula above if there
were a way of determining what value of $i$ gives the largest
answer.
\m {\bf 24. Remark.} Since for all $j$ there is an $i$ so that the pushforwards to $Y_0$  of ${\cal O}_U(jK_U)$
is dual to the pullback  from $M$ to $Y_0$  of ${\cal F}_{(n+2)^i-1}^{\otimes j}$ times
the inverse of the pullback modulo torsion of ${\cal F}_{(n+2)^i}^{\otimes j}$
then it is a reflexive sheaf for all  $j$ on any bounded open subset of $Y_0.$
According to a question of  [7] then  it is not known whether  $Y_0$ has canonical singularities.

\vfill\eject \centerline{12. Appendix}
\m\m We will finish this note by describing the
 proof of 
 properness of $V\to M.$
In this section we assume $M$ to be normal and locally algebraic. Refer to [21,22,23]  for the full list of references and full statements
of theorems. Recall $\pi:M'\to M$ is a resolution.
\m Kempf describes things this way in his article which happens to
be in the Springer
Lecture notes [6] about toroidal embeddings: although we have not assumed
$0=R^i\pi_*\Lambda^{n}\Omega_{ M'}$ for $i\ge 1,$
Grauert-Riemenshneider vanishing says this is true.
\m We may assume $M$ is a closed analytic subvariety
of a disk  $A$
by an embedding $i$
of codimension $c.$ The sheaf $\Lambda^{dim(A)}\Omega_A$ is isomorphic
to ${\cal O}_A$ but not canonically and we shall be needlessly
rigorous about the notation by distinguishing them.
 Since all higher derived functors
vanish,  duality simply implies
${\cal E}xt^c_{{\cal O}_A}(-, \Lambda^{dim(A)}\Omega_A)$
interchange $i_*\pi_*{\cal O}_{M' }$ and
$i_*\pi_*\Lambda^{n}\Omega_{M'}.$
The isomorphism ${\cal O}_M\to \pi_*{\cal O}_{M'}$ coming
from normality of $M$ gives when we apply $i_*$ and apply our ${\cal E}xt$
functor
 $$i_*\pi_*\Lambda^{n}
\Omega_{M'}\cong {\cal E}xt^c_{{\cal O}_A}(i_*\pi_*{\cal O}_{M'}, \Lambda^{dim(A)
}\Omega_A)$$
$$=  {\cal E}xt^c(i_*{\cal O}_M, \Lambda^{dim(A)}\Omega_A)$$
$$=i_*{\cal E}xt^c({\cal O}_M, \Lambda^{dim(A)}\Omega_A)$$
Removing the $i_*$ this is the double dual of the highest exterior
power of differentials of $M, $  ie the result of removing all but
codimension one primary components.
Thus this is the pushforward of the highest exterior
power of the differentials of $M'$
$$\cong i_*\pi_*\Lambda^{n}\Omega_{M'}.$$
\m The definition of rational singularities
doesn't imply the $\pi_*{\cal O}_{M'}(iK_{M'})$ for $i\ge 2$ are the sheaves associated
to the divisors $iK$ where $K$ is the canonical divisor. 
 The duality argument above does
not extend  to the higher values of $i.$  The condition that
$\pi_*{\cal O}_{M'}(iK_{M'})={\cal O}_V(iK_M),$ with the extra
condition of ${\mathbb Q}$-Cartier (but note it is
stated in [7] that the extra condition may be
automatic) is the definition of $M$ having canonical singularities.
\m  Always one inclusion holds
 $\pi_*{\cal O}_{M'}(iK_{M'}) \subset {\cal O}_V(iK_M),$ and this is natural
independent of choice of $M'.$  This is because $M$ is normal,
and any Weil divisor is then Cartier except on a locus
of codimension at least two. To see this, we can assume
we are talking about an irreducible Weil divisor which for
a normal variety corresponds to a symboic power of a height one
prime ideal. Choosing an element of the local ring at
that prime which generates the corresponding  power
of the prime ideal, this defines the same divisor except on
a locus of vanishing of a single element, and the divisor
of this element in each affine coordinate chart of a finite
cover of the divisor  meets the divisor
in codimension at least two.
\m  Then the blowup of any sheaf of ideals, by taking a primary
decomposition, we see it is an isomorphism away from a codimension
two locus. Or we could have chosen our resolution to be an isomorphism
away from the codimension at least two singular locus. In any
case then the direct image of ${\cal O}_{ M'}(iK_{M'})$
agrees with
${\cal O}_V(iK_M)$ except in codimension at last two. The latter
has no primary components  with associated primes of height other than one,
and so cannot be made larger without affecting something codimension one.
\m Shepherd Barron in dimension 3 [8] and  Elkik
in dimension $\ge 4$ [9] answered one of the questions in
[7] namely that the Cohen Macaulay property, vanishing
of all but one ${\cal E}xt$, does follow from
the conditions of canonical singularities.
\m
 Assume the relative minimal model
program  as described in the
conjecture of  [11]
holds for the resolution of $\pi:M'\to M.$ 
So we are assuming  Mori's sequence of contractions
(followed by passing to the relative canonical model) [13] leading from $M'$ to a variety $q:M_t
\to M$
 which has $K_{M_t}$ relatively nef over $M$ and terminal
singularities. The relative basepoint free theorem as described in the textbook [21], page 94, even in the analytic case, gives some $iK_{M_t}$
relatively basepoint free. 
Choose a resolution
$h:M''\to M_t$ of $M_t$; one has
$$h_*{\cal O}(iK_{M''})={\cal O}(iK_{M_t})$$
By relative basepoint freeness this is
$$=q^*q_*{\cal O}_{M_t}(iK_{M_t})$$
$$=q^*q_*h_*{\cal O}(iK_{M''})$$
By uniqueness (which uses normality of $M,$ that $iK_{M'}$ is Cartier
and contravariant) this is
$$=q^*\pi_*{\cal O}(iK_{ M'}).$$
$$=q^*\omega_i$$
It follows from these (inefficient) formulas
that $$q_*{\cal O}_{M_t}(iK_{M_t})=\omega_i.$$
\m The relative basepoint freeness holds with $i$ replaced
by $2i, 3i, ...$ because all it is saying is (assuming as
we may that $M$ is affine) that at each point of $M_t$ there
is a global section with the correct order of pole there,
and powers of global sections give the correct order
for multiples of the divisor. Therefore there is a continuous function
$M_t\to V$ over $M$ and  compact subsets of $M$ lift to compact
subsets of $M_t$ which remain compact in $V.$

\vfill\eject \centerline{
 References}
\m 

\noindent
1. Michael Atiyah, Complex analytic connections in fibre Bundles, Transactions of the American Mathematical Society, Vol. 85, No. 1,  181-207  (1957)
\vspace{6pt}

\noindent
2. Alexandre Grothendieck, Elements de geometrie  algebrique II, 
publ IHES 8, 5-222 (1961)
\vspace{6pt}

\noindent
3. Heisuke Hironaka, Resolution of singularities of an algebraic variety over a field of characteristic zero II, Annals of Mathematics 79, 205-326 (1964)
\vspace{6pt}

\noindent
4. John  Nash, unpublished (circa 1965)
\vspace{6pt}

\noindent
5. Michel Artin, Algebraic approximation of structures over
complete local rings, IHES 36, 23-58 (1969)
\vspace{6pt}

\noindent 
6. George Kempf , Finn Knudsen, David Mumford, Bernard Saint-Donat,
 Toroidal Embeddings 1,  Lecture Notes in Mathematics , Vol. 339
(1973)
\vspace{6pt}

\noindent
7. Miles Reid, Canonical 3-folds, proceedings of the Angiers `Journees de Geometrie Algebrique' (1979)
\vspace{6pt}
 
\noindent 
8. Nick Shepherd-Barron, PhD thesis, Warwick University.
\vspace{6pt}

\noindent 
9. Renee Elkik   Rationalite des singularites canoniques, Invent. Math. 64 (1981) 
\vspace{6pt}

\noindent
10. Heisuke Hironaka, On Nash blowing-up, in Arithmetic and geometry, Vol. II, 103.111, Progr. Math., 36, Birkhauser, Boston, Mass., (1983) 
\vspace{6pt}

\noindent
11. Miles Reid,  Decomposition of Toric Morphisms, 
in Arithmetic and Geometry, vol. II: Geometry, Prog. Math. 36,. 
395-418 (1983)
\vspace{6pt}

\noindent
12. Helene Esnault, Eckart Viehweg, Logarithmic deRham complexes and
vanishing theorems, Invent. Math 86, 191-194 (1986)
\vspace{6pt}

\noindent
13. Shigefumi Mori, Birational classification of algebraic threefolds,
International Congress of Mathematicians, Kyoto (1990)
\vfill\eject\ 

\noindent
14. Mark Spivakovsky, Sandwiched surface singularities and desingularization
of surfaces by normalized Nash transformations,
 Annals of Mathematics,  Vol. 131, No. 3, 411-491 (1993)
\vspace{6pt}

\noindent
15. John Nash, Arc structure of singularities, Duke Math. J. Volume 81, 
Number 1, 31-38  (1995)
\vspace{6pt}

\noindent
16. Mauro Beltrametti, Andrew Sommese, The adjunction theory of complex 
projective varieties, DeGruyter (1995)
\vspace{6pt}

\noindent
17. John Moody, On resolving singularities,  Journal of the London Mathematical Society, 64 , pp 548-564 (2001)
\vspace{6pt}

\noindent
18. Shihoko Ishii, Janos Kollar, The Nash problem on arc families of 
singularities Duke Math. J. Volume 120, Number 3,601-620 (2003)
\vspace{6pt}

\noindent
19. John Moody, Finite generation and the Gauss Process, preprint (2004)
\vspace{6pt}

\noindent
20. Janos Kollar, Lectures on resolution of singularities, 
Annals of Mathematics Studies 166 (2007)
\vspace{6pt}

\noindent 
21. Janos Kollar, Shigefumi Mori, Birational geometry of
algebraic varieties,
Cambridge University Press (United Kingdom), (2008) 
\vspace{6pt}

\noindent
22. James McKernan, Recent advances in the MMP after Shokurov, II,
(lecture)
\vspace{6pt}

\noindent 
23. Caucher Birkar,  Paolo Cascini, Christopher Hacon, James McKernan,
Existence of minimal models for varieties of log general types, J Amer Math Soc  23 (2010) 405-468
\end{document}